%% file: agt-1-24.tex
\documentclass{gtart}
\input agtout

\lognumber{24}
\volumenumber{1}
\volumeyear{2001}
\papernumber{24}
\pagenumbers{469}{489}
\received{7 August 2001}
\accepted{6 September 2001}
\published{9 September 2001}

\usepackage{amssymb,amscd,amsmath}

\newtheorem{thm}{Theorem}[section]
\newtheorem{lem}[thm]{Lemma}
\newtheorem{prop}[thm]{Proposition}
\newtheorem{cor}[thm]{Corollary}
\theoremstyle{definition}
\newtheorem{defn}[thm]{Definition}
\newtheorem*{rem}{Remark}

\newcommand{\tam}{\mathbb Z}
\newcommand{\kar}{\mathbb C}
\newcommand{\Real}{\mathbb R}
\newcommand{\Euc}{\mathbb E}
\newcommand{\Hyp}{\mathbb H}

\newcommand{\To}{\longrightarrow}

\newcommand{\G}{\pi_1(M)}
\newcommand{\product}{S^1 \times M}

\begin{document}

\title{Lefschetz fibrations, complex structures and\\Seifert 
fibrations on $\product^3$}
\authors{Tolga Etg\"{u}}
\covertitle{Lefschetz fibrations, complex structures and\\Seifert 
fibrations on $S^1\times M^3$}
\coverauthors{Tolga Etg\noexpand\"{u}}
\asciititle{Lefschetz fibrations, complex structures and Seifert 
fibrations on S^1 X M^3}
\asciiauthors{Tolga Etgu}
\address{Department of Mathematics\\University of California at 
Berkeley\\Berkeley, CA 94720, USA}
\email{tolga@math.berkeley.edu}

\begin{abstract}

We consider product 4--manifolds $\product$, where $M$ is a closed, connected and
oriented 3--manifold.  We prove that if $\product$ admits a complex structure or a
Lefschetz or Seifert fibration, then the following statement is true:
\smallskip

\cl{\sl{$\product$ admits a symplectic structure if and only if $M$ fibers
 over $S^1$,}}

under the additional assumption that $M$ has no fake 3--cells. We also discuss the
relationship between the geometry of $M$ and complex structures and Seifert fibrations
on $\product$.

\end{abstract}

\asciiabstract{We consider product 4--manifolds S^1 X M, where M is a
closed, connected and oriented 3-manifold.  We prove that if S^1 X M
admits a complex structure or a Lefschetz or Seifert fibration, then
the following statement is true: S^1 X M admits a symplectic structure
if and only if fibers over S^1, under the additional assumption that
M has no fake 3-cells. We also discuss the relationship between the
geometry of M and complex structures and Seifert fibrations on
S^1 X M.}

\primaryclass{57M50, 57R17, 57R57} \secondaryclass{53C15, 32Q55}
\keywords{Product 4--manifold, Lefschetz fibration, symplectic
manifold, Seiberg--Witten invariant, complex surface, Seifert
fibration}
\asciikeywords{Product 4-manifold, Lefschetz fibration, symplectic
manifold, Seiberg-Witten invariant, complex surface, Seifert
fibration}

\maketitle

\section{Introduction}

A closed, oriented, smooth 4--manifold $X$ which fibers over a Riemann surface admits
a symplectic structure unless the fiber class is torsion in $H_2(X; \tam)$. In
particular, a fibration of a closed, oriented 3--manifold $M$ over $S^1$ induces a
symplectic form on $\product$.

\medskip
\noindent {\bf{Conjecture T}}\qua {\sl{Let $M$ be a closed, oriented 3--manifold such
that $\product$ admits a symplectic structure. Then $M$ fibers over $S^1$.}}
\medskip

This conjecture was first stated by Taubes \cite{tau} and is still open. Recent work
of Chen and Matveyev \cite{chemat} proves that it holds when $M$ has no fake 3--cells,
$\product$ admits a symplectic structure and a Lefschetz fibration with symplectic
fibers.

In this paper, we generalize Chen and Matveyev's result proving that the conjecture
holds when $\product$ admits an arbitrary Lefschetz fibration (possibly with
nonsymplectic fibers). More generally, we prove the following:

\begin{thm} \label{mainthm}

Suppose $M$ is a closed 3--manifold without a fake 3--cell.

{\rm(L)}\qua If $\product$ admits a Lefschetz fibration, then Conjecture T holds.

{\rm(S)}\qua If $\product$ admits a Seifert fibration, then Conjecture T holds.

{\rm(K)}\qua If $\product$ admits a K\"{a}hler structure, then Conjecture T holds.

{\rm(C)}\qua If $\product$ admits a complex structure, then Conjecture T holds.

\end{thm}
Here, a fake 3--cell means a compact, contractible 3--manifold which is not
homeomorphic to $D^3$. Note that the Poincar\'{e} conjecture implies that there is no
fake 3--cell.

\begin{rem}
We'll see that a nonsymplectic Lefschetz fibration on a product 4--manifold has no
singular fibers and has fiber a torus. Since a Seifert fibration can be thought of as
a $T^2$--fibration with multiple fibers, (S) is a further generalization of (L).  
Statement (C)
is clearly a generalization of (K). Note that all (symplectic) product 4--manifolds
which admit complex structures turn out to be Seifert fibered. This means that all
other statements follow from (S) using the result of Chen and Matveyev on symplectic 
Lefschetz fibrations.
\end{rem}

In the remark above and the rest of the paper, by a product 4--manifold we mean the
product of $S^1$ with a (compact, oriented, connected) 3--manifold.

In order to prove Theorem \ref{mainthm}, besides other techniques, we use
classification results on complex surfaces and Lefschetz fibered 4--manifolds and
apply them to product manifolds. In particular, we get results on the classification
of product 4--manifolds which admit certain structures or fibrations and interesting
relations between the geometry of $M$ and complex structures and Seifert fibrations on
$\product$.

\begin{rem}
In their paper \cite{geigo} on taut contact circles on 3--manifolds, Geiges and
Gonzalo classified product 4--manifolds carrying complex structures with respect to
which the obvious circle action is holomorphic. Since we don't require this action to
be holomorphic and we are mainly interested in the symplectic structure on product
manifolds, we prove different type of results even though we use similar methods.
\end{rem}

\begin{rem}
As a consequence of Theorem \ref{mainthm} we see that when $M$ is a nonhyperbolic
geometric 3--manifold Conjecture T holds. On the other hand, assuming Thurston's
conjecture on the geometrization of 3--manifolds, if $\product$ admits a symplectic
structure, then $M$ is prime (see \cite{mccar} or \cite{vid}). So it might be
interesting to try to prove Conjecture T (at least up to the geometrization
conjecture) by first proving it when $M$ is hyperbolic, then considering geometric
3--manifolds with boundary (disjoint union of tori) and finally using Seiberg--Witten
theory of 4--manifolds glued along $T^3$.
\end{rem}
In the next section we recall definitions and some basic theorems on Lefschetz
fibrations, complex surfaces, Seiberg--Witten invariants, Seifert fibrations and
geometric structures on 3-- and 4--manifolds. In Section 3, we discuss nonsymplectic
Lefschetz fibrations on $\product$. By using the Seiberg--Witten theory of symplectic
4--manifolds and $S^1$--bundles over surfaces, we prove (L) of Theorem \ref{mainthm}
in Section 4. In Section 5, product 4--manifolds which admit complex structures are
considered and (K) is proved first. As a result of a slightly more careful
investigation we prove (C). Finally we consider Seifert fibered 4--manifolds and prove
(S). In the last section, we discuss the relation between various structures and
fibrations on $\product$ and $M$.

In this paper, by a fiber bundle we mean a locally trivial one and an $F$--bundle
means a (locally trivial) fiber bundle with fiber $F$. All fibrations (of any kind)
are oriented and all manifolds are compact, smooth, oriented and connected, unless
stated otherwise.

\noindent {\bf{Acknowledgment}}\qua The author is grateful to his thesis advisor Rob
Kirby for numerous discussions.

\section{Background}

Let us first state some topological information on $\product$.
\begin{lem} \label{basics}
Let $M$ be a closed, oriented and connected 3--manifold. Then $X = \product$ is a spin
4--manifold with $\sigma (X)=\chi (X) = 0$, $b_{\pm}(X) = b_1(M)$ (in particular,
$b_2(X)$ is even), where $\sigma$, $\chi$ and $b_*$ denote the signature, Euler
characteristic and the corresponding Betti number, respectively.
\end{lem}
\begin{proof}
Both $S^1$ and $M$ are spin, so $X$ is spin. Since $\chi (S^1) =0 $,
the Euler characteristic of $X$ vanishes. The boundary of $D^2 \times M$ is $X$, so $\sigma (X) = 0$.
The facts about the Betti numbers follow easily from the definitions of $\sigma$ and $\chi$ in terms of Betti
numbers.
\end{proof}

\subsection{Lefschetz fibrations and pencils}

\begin{defn}
A Lefschetz fibration on a compact, connected, oriented and smooth 4--manifold $X$ is
a smooth map $\pi \co X \To \Sigma $, where $\Sigma$ is a compact, connected, oriented
surface and $\pi^{-1} (\partial \Sigma ) = \partial X $, such that each critical point
of $\pi$ lies in the interior of $X$ and has an orientation-preserving coordinate
chart on which $\pi (z_1, z_2) = z_1^2+z_2^2$ relative to a suitable smooth chart on
$\Sigma$.
\end{defn}

\begin{defn}
A Lefschetz pencil on a closed, connected, oriented, smooth 4--manifold $X$ is a
non-empty finite subset $B$ of $X$ called the base locus, together with a smooth map
$\pi \co X - B \To \kar P^1$ such that each point $b \in B$ has an
orientation-preserving coordinate chart in which $\pi$ is given by the
projectivization $\kar^2 - \{ 0 \} \To \kar P^1$, and each critical point has a local
coordinate chart as in the definition of a Lefschetz fibration above.
\end{defn}

\begin{defn}
A Lefschetz fibration is called relatively minimal if no fiber contains an exceptional
sphere, in other words it cannot be obtained by blowing up another Lefschetz
fibration.
\end{defn}

\begin{defn}
A Lefschetz fibration is called a symplectic Lefschetz fibration if the total space
admits a symplectic structure such that generic fibers are symplectic submanifolds,
otherwise it is called nonsymplectic.
\end{defn}

\begin{thm}[Gompf] \label{symplf}
A Lefschetz fibration on a 4--manifold $X$ is symplectic if and only if the homologous
class of the fiber is not torsion in $H_2(X;\tam )$.
\end{thm}

The close relation between Lefschetz fibrations and symplectic structures is stated in
the following theorems.

\begin{thm}[Donaldson \cite{don}] Every symplectic
4--manifold admits a Lefschetz pencil by symplectic surfaces.
\end{thm}

\begin{thm}[Gompf \cite{gomsti}] If a 4--manifold admits
a Lefschetz pencil (with non-empty base locus), then it admits a symplectic structure.
\end{thm}

It is necessary that the base locus is non-empty as we have examples of 4--manifolds,
e.g.\ $S^1 \times S^3$, which admit Lefschetz fibrations over $S^2$ but no symplectic
structure.

If a manifold admits a Lefschetz pencil, then one can blow-up the points of the base
locus and construct a Lefschetz fibration (over $S^2$). So Donaldson's theorem implies
that every symplectic 4--manifold has a blow-up which admits a Lefschetz fibration.
Even though it is always possible to put a Lefschetz pencil on a symplectic $\product$
it may not be possible to find a Lefschetz fibration on it. Note that a blow-up of
$\product$ can never be a product.

For more details on Lefschetz pencils and fibrations see \cite{gomsti}.

\subsection{Seiberg--Witten invariants}

Let $X$ be a closed, oriented, connected and homology oriented 4--manifold with
$b_+(X)
> 0$. The Seiberg--Witten invariant $SW$ of a $Spin_c$ structure on $X$ was first extracted
from monopole equations by Witten in \cite{wit}. If $b_+(X) > 1$, then $SW$ is an
integer-valued diffeomorphism invariant of $X$. When $b_+(X)=1$ it may depend on the
chosen metric. The Seiberg--Witten invariant of a $Spin_c$ structure $\xi$ on $X$ is
denoted by $SW_X(\xi)$. We call $\alpha \in H^2(X; \tam)$ a basic class if there
exists a $Spin_c$ structure $\xi$ such that $SW_X(\xi) \neq 0 $ with $c_1(\det (\xi))
= \alpha$, where $\det (\xi)$ denotes the determinant (complex) line bundle of $\xi$.
If there is no 2--torsion in $H^2(X; \tam)$, then there is a unique $Spin_c$ structure
$\xi$ with $c_1(\det (\xi)) = \alpha$ for any characteristic class $\alpha \in H^2(X;
\tam)$. In general, the set of isomorphism classes of $Spin_c$ structures on $X$ is an
affine space modeled on $H^2(X; \tam)$.

Seiberg--Witten invariants of 3--dimensional manifolds are defined similarly. As we
state in Section 4, Seiberg--Witten invariants of a 3--manifold $M$ carry exactly the
same information as those of $\product$ at least when $b_1(M)
> 1$. The reader is referred to \cite{kron} and \cite{okotel}
for the theory of Seiberg--Witten invariants in dimension 3.

We have the following important theorem on the Seiberg--Witten invariants of
symplectic manifolds.

\begin{thm}[Taubes \cite{tau2}, \cite{tau3}]  \label{swsymp}
Let $X$ be a closed 4--manifold with $b_+
> 1$ and a symplectic form $\omega$. Then there is a canonical $Spin_c$ structure $\xi$ on
$X$ such that $SW_X(\xi) = \pm 1 $ and $\det(\xi)$ is the canonical line bundle $K$ of
$(X, \omega)$.

Moreover,
$$ 0 \leq | \alpha \cdot [ \omega ] | \leq | c_1 (K) \cdot [ \omega ] | \ , $$
where $\alpha$ is any basic class; $0= \alpha \cdot [ \omega ] $ if and only if
$\alpha =0$; $| \alpha \cdot [ \omega ] | = | c_1 (K) \cdot [ \omega ] |$ if and only
if $\alpha = \pm c_1(K) $.
\end{thm}

See \cite{gomsti},  \cite{mor} and \cite{kron} for more details on Seiberg--Witten
invariants of 4-manifolds.

\subsection{Geometric structures and the geometrization conjecture}
\begin{defn}
A metric on a manifold is called locally homogeneous if any pair of points can be
mapped to each other by isometries of open neighborhoods.
\end{defn}

\begin{defn}
A manifold is called geometric if it admits a complete, locally homogeneous metric.
\end{defn}

\begin{defn}
A simply connected geometric manifold together with the isometry group corresponding
to a complete (locally) homogeneous metric is called a geometry.
\end{defn}

Up to isometry, there are eight 3--dimensional and nineteen 4--dimensional geometries
with compact quotients. These are classified by Thurston and Filipkiewicz \cite{fil}
respectively. See \cite{sco} and \cite{wall} for detailed discussions on 3-- and
4--dimensional geometries.

A manifold is called prime if it cannot be written as the connected sum of two
manifolds none of which is a sphere. In \cite{mil} Milnor showed that, up to
homeomorphism and the permutation of the summands, there is a unique way to write a
compact, oriented 3--manifold as the connected sum of prime manifolds. There is also a
reasonably canonical way to cut compact, prime 3--manifolds along tori into pieces
which no longer have embedded tori in them other than their boundary components (up to
homology). Thurston's geometrization conjecture asserts that these pieces should all
be geometric.

\subsection{Seifert fibered spaces}

A trivial fibered solid torus is $S^1 \times D^2$ with the product foliation by
circles. A fibered solid torus is a solid torus with a foliation by circles that is
finitely covered by a trivial fibered solid torus. It can be constructed by gluing two
ends $D^2 \times \{ 0 \} $ and $D^2 \times \{ 1 \}$ of $D^2 \times I$ after a $q/p$
rotation.

A Seifert fibered space is a 3--manifold with a decomposition into disjoint circles
such that each circle has a neighborhood isomorphic to a fibered solid torus. A circle
bundle over a surface is naturally a Seifert fibered space. By identifying each of
these circles with a point, we can consider a Seifert fibered space as a fibration
over a 2--orbifold base. Such a fibration is called a Seifert fibration. Fibers of a
Seifert fibration are obviously circles and singularities of the base orbifold
correspond to the fibers without trivial fibered solid torus neighborhoods. A fiber is
called regular if it projects to a nonsingular point of the base, otherwise it is
called a multiple fiber.

\begin{lem}[cf.\ Lemma 3.2 in \cite{sco}] \label{shs}
Suppose $M$ admits a Seifert fibration over a 2--orbifold $X$. Then there is a short
exact sequence
$$ 1 \To G \To \pi_1(M) \To \pi^{orb}_1(X) \To 1 \ ,$$
where $G$ denotes the cyclic subgroup of $\G$ generated by a regular fiber and
$\pi^{orb}_1(X)$ denotes the fundamental group of $X$ as an orbifold. The subgroup $G$ is
infinite
except in cases where $M$ is covered by $S^3$.
\end{lem}

Note that a presentation for $\pi^{orb}_1(X)$ is
$$ \big\langle a_1 , b_1 , \dots , a_g , b_g , x_1 , \dots , x_n\ \big|\ x_i^{p_i} = 1 ,
\prod^g_{i=1} [ a_i , b_i ] \cdot \prod^n_{i=1}x_i = 1 \big\rangle \ ,$$ 
where $g$ is the genus of the underlying surface of $X$,
assuming $X$ is closed and orientable with $n$ singular points of multiplicities $p_1,
\dots p_n$. The Euler characteristic $\chi(X)$ of such a 2--orbifold $X$ is defined by 
$$ \chi (X) =2-2g - \sum^n_{i=1} \left( 1 - \frac{1}{p_i} \right) \ .$$ 
An orbifold is called spherical (Euclidean or hyperbolic) if its Euler characteristic
is positive (zero ornegative).

For more details on Seifert fibered spaces see \cite{orl} and \cite{neuray}. For
geometric structures on Seifert fibered spaces see \cite{sco}.

\subsection{Seifert fibered 4--manifolds}

A Seifert fibration on a 4--manifold is analogous to a Seifert fibration on a
3--manifold.

\begin{defn}
A smooth map $\pi \co X \To \Sigma$ from a smooth 4--manifold $X$ to a surface
$\Sigma$ is called a Seifert fibration if there exists a finite set of isolated points
$B$ in $\Sigma$ such that the restriction of $\pi$ to $\pi^{-1} (\Sigma -B)$ is a
torus bundle and for each element $b \in B$, $\pi^{-1}(b)$ has a tubular neighborhood
diffeomorphic to the product of a fibered solid torus with a circle.
\end{defn}

A Seifert fibration can be thought of as a torus fibration over a 2--orbifold. In the
complex category it corresponds to an elliptic fibration without singular fibers
(possibly with multiple ones). If a 4--manifold admits a Seifert fibration it is
called a Seifert 4--manifold. We have analogous statements for Seifert fibered
4--manifolds to most of the properties of Seifert fibered spaces, e.g.\ Lemma
\ref{shs}. See \cite{wall} and \cite{wall2} for geometric structures on elliptic
surfaces without singular fibers, \cite{ue1} and \cite{ue2} for a general picture of
Seifert 4--manifolds in terms of geometric structures.

\section{Nonsymplectic Lefschetz fibrations on $\product$}

In this section our aim is to show that nonsymplectic Lefschetz fibrations on
$\product$ are in fact locally trivial torus bundles. We also investigate which of
these fibrations have symplectic total spaces and which of them give rise to
fibrations of $M$ over $S^1$.

\begin{thm}[Chen-Matveyev \cite{chemat}] \label{slfimpfib}
Let $\pi$ be a symplectic Lefschetz fibration on $\product$, where $M$ is a closed,
connected, oriented 3--manifold without any fake 3--cells. Then there exists a
fibration $p$ on $M$ over $S^1$. Moreover, the symplectic structure with which $\pi$
is compatible is deformation equivalent (up to self-diffeomorphisms of $\product$) to
the canonical symplectic structure associated to the fibration $Id \times p  \co
\product \rightarrow S^1 \times S^1$.
\end{thm}

The symplectic form (canonical up to deformation equivalence) on the total space of a
surface bundle over a compact, oriented surface is obtained by extending a symplectic
form on a fiber and adding a (sufficiently large) multiple of the pullback of a
symplectic form on the base to it (see \cite{thu} and \cite{mcdsal} for details and
more general cases). The following lemma plays a crucial role in the proof of the
theorem above.
\begin{lem}{\rm\cite{chemat}}\qua \label{lemvancyc}
Let $\pi$ be a symplectic Lefschetz fibration on $\product$, where $M$ is a closed,
connected, oriented 3--manifold. Then $\pi$ doesn't have any critical points.
\end{lem}

First of all, we give the following generalization of this lemma.

\begin{lem}  \label{nsltb}
Let $\pi$ be a Lefschetz fibration on $\product$, where $M$ is a closed, connected,
oriented 3--manifold. Then $\pi$ is a fiber bundle. If $\pi$ is not symplectic, then
it is a torus bundle.
\end{lem}

\begin{proof}
We only need to consider the case where $\pi$ is not symplectic, i.e. fibers are not
symplectic submanifolds of $X= \product$. By Theorem \ref{symplf} the fiber class $[ F
]$ is torsion in $ H_2 ( X; \tam ) $. This is possible only if $F$ is a torus since
otherwise
$$ 0 \neq \chi (F) = \langle e (T F), [ F ] \rangle \ .$$
Note that $e (T F)$ extends to $H^2(X ; \tam)$ since $TF$ is the pull-back (by the
inclusion $F \hookrightarrow X$) of the vertical (with respect to $\pi$) subbundle of
$T X$. On the other hand, the Euler characteristic of the total space of a Lefschetz
fibration is equal to the product of the those of the base and the fiber plus the
number of vanishing cycles (assuming there is a unique singular point on each fiber).
In our case this leads to
$$ 0= \chi (\product ) = \chi (T^2) \cdot \chi (B) + \# \{ \hbox{vanishing cycles}\} \ .$$
Hence there are no vanishing cycles. Therefore $\pi$ is a torus bundle.
\end{proof}

This lemma shows that nonsymplectic Lefschetz fibrations on $\product$ are all torus
bundles over Riemann surfaces. We investigate these bundles in three groups according
to the genera of their bases.

\begin{lem} \label{t2overs2}
Let  $\product$ be the total space of a nontrivial $T^2$--bundle over $S^2$. Then
$\product$ carries no symplectic form.
\end{lem}

\begin{proof}
Since the torus bundle is nontrivial, $b_1(\product) < 2$ and therefore 
$b_2 (\product) = 2 \cdot b_1(M) = 0$. Hence all closed 2--forms on $\product$ are degenerate.
\end{proof}

\begin{rem}
As we mentioned before, a fibration of $M$ over $S^1$ induces a symplectic form on
$\product$. Therefore, when $\product$ is as in the lemma $M$ doesn't fiber over the
circle.
\end{rem}

We have a totally different picture for $T^2$--bundles over $T^2$.

\begin{thm}[Geiges \cite{gei}] \label{t2overt2}
Let $X$ be the total space of an oriented $T^2$--bundle over $T^2$. Then $X$ admits a
symplectic structure. Moreover, there exists a symplectic $T^2$--bundle over $T^2$
with total space $X$ unless $X$ is the total space of a nontrivial $S^1$--bundle over
the total space of a nontrivial $S^1$--bundle over $T^2$.
\end{thm}

Let $X$ be an exception, i.e. a twisted circle bundle over a twisted circle bundle
over the torus. Then $b_1(X) = b_2(X) = 2 $. Moreover, $H^1_{DR} (X; \Real)$ is generated
by $[ \alpha ]$ and $[ \beta ] $, where $\alpha $ and $\beta$ are closed 1--forms on
$X$ such that $n \cdot \alpha \wedge \beta = d \gamma$, where $n$ is the Euler number of the
(nontrivial) $S^1$--bundle over $T^2$ and $\gamma$ is a 1--form on $X$ (see \cite{fgg}
for details). In particular, $(H^1(X; \Real))^{\cup 2} = 0$, where $(H^1(X;
\Real))^{\cup 2}$ denotes the image of the cup product of $H^1(X; \Real)$ with itself.
On the other hand, $H^1(\product ; \Real) \cong H^1(S^1; \Real) \oplus H^1(M; \Real )
$ and obviously $(H^1(\product; \Real))^{\cup 2} \neq 0$. Therefore we have the
following corollary.

\begin{cor} \label{t2ovt2slf}
If $\product$ is the total space of a $T^2$--bundle over $T^2$, then $\product$ admits
a symplectic Lefschetz fibration.
\end{cor}

For $T^2$--bundles over higher genus surfaces we have

\begin{lem} \label{sl2dnfib}
Let  $\product$ be the total space of a $T^2$--bundle  over $B$, where $B$ is a
closed, oriented surface of genus $\geq 2$. Also assume that $M$ has no fake 3--cells.
Then $M$ fibers over the circle if and only if the torus bundle is trivial.
\end{lem}

We are going to use the following lemma to prove the one above.

\begin{lem}[cf.\ \cite{orl} Theorem 7.2.4] \label{lemsl2dnfib}
Let $M$ be a closed,oriented 3--manifold which is the total space of a circle bundle
over a closed, oriented surface $B$ of genus $\geq 2$. Then $M$ fibers over the circle if
and only if $M = S^1 \times B$.
\end{lem}

\begin{proof}
Recall that $\pi_1(M)$ has the presentation
$$\big\langle a_1 , b_1 , \ldots , a_g, b_g ,
\alpha\ \big|\ [a_i , \alpha ] = [b_i , \alpha ] = 1 , [ a_1,b_1 ] \cdots [a_g , b_g] =
\alpha^k \big\rangle \ ,$$ 
where $g= genus (B)$ and $k$ is the Euler number of the
$S^1$--bundle. In particular, $ H_1 (M) \cong \tam^{2g+1} $ if $k=0$ and $H_1(M) \cong
\tam^{2g} \oplus \tam_{|k|}$ otherwise.

We also have the following commutative diagram of exact sequences
\[ \begin{CD}
0 @>>> \pi_1(S^1) @>j_{\#}>>  \pi_1(M) @>>> \pi_1 (B) @>>> 1 \\
 @. @VV\cong V @VVV @VVV @.  \\
 @. H_1(S^1) @>j_{*}>> H_1 (M) @>>> H_1(B) @>>> 0
\end{CD} \]
where vertical maps are Hurewicz epimorphisms. Note that the homomorphism $j_* $ is injective 
if and only if $Im(j_{\#}) \cap [\G:\G] = \{ 1 \}$. Now suppose that $F \To M \To S^1$ is a
fibration.
There exists a normal subgroup $N \cong \pi_1(F) $ in $\G$ such that
$\G /N \cong
\tam$. Assume that there exists an element $u \in N \backslash \{ 1 \} $ such
that
$u=j_{\#} (v)$. Then there is a normal infinite cyclic subgroup (generated by $u$) in
$N$ and this implies that $F$ is a torus, but $M$ cannot be the total space of a torus
bundle over the circle since $b_1(M) \geq 2g \geq 4 > 3 $. Therefore $Im(j_{\#}) \cap
N = \{ 1 \}$. On the other hand, $[\G:\G] \subset N$ because $\G /N \cong \tam$. So
$Im(j_{\#}) \cap [\G:\G] = \{ 1 \}$, $j_*$ is injective and we have the short exact
sequence
$$ 0 \To H_1(S^1) \To H_1 (M) \To H_1(B) \To 0 $$
which clearly splits. Hence $b_1(M)=2g+1$ and $M$ is the product $S^1 \times B$.
\end{proof}

\proof[Proof of Lemma \ref{sl2dnfib}] We have the homotopy sequence of the
$T^2$--bundle
\begin{equation} \label{t2homseq}
 0 \To  \pi_1(T^2) \xrightarrow{j_{\#}} \pi_1(\product )
\xrightarrow{\pi_{\#}} \pi_1(B)\To 1\ .
\end{equation}
Let $u$ be a generator of $\pi_1(S^1 \times {pt} ) $. Assume that $\pi_{\#} (u) = v
\neq 1 \in  \pi_1(B) $. Then $v$ generates a normal cyclic subgroup in $\pi_1(B)$ and
this contradicts the fact that $genus (B) \geq 2$. Therefore $u \in ker(\pi_{\#}) = im
(j_{\#}) $, where $j$ is the inclusion map. Let $a$ be $j_{\#}^{-1} (u) $. We can
find another element $b \in \pi_1(T^2) $ such that the restriction of $j_{\#}$ to the
subgroup $\langle b \rangle$ generated by $b$ gives the short exact sequence

\begin{equation} \label{s1homseq}
0 \To \langle b \rangle \To  \pi_1(M) \To \pi_1 (B) \To 1\ .
\end{equation}

By Theorem 11.10 in \cite{hem} $M$ admits an $S^1$--bundle over $B$ (we use the
assumption that $M$ has no fake 3--cells). Lemma \ref{lemsl2dnfib} finishes the proof.
\endproof

We should note that the idea of extracting (\ref{s1homseq}) from (\ref{t2homseq}) was
first used in \cite{chemat}.

\begin{prop} \label{oldmainthm}
Suppose $\product$ admits a nonsymplectic Lefschetz fibration, where $M$ is a closed,
oriented 3--manifold. If the base space of the fibration is a torus, then $\product$
admits a symplectic form and a symplectic Lefschetz fibration. Otherwise $M$ doesn't
fiber over $S^1$ or it has a fake 3--cell.
\end{prop}

\begin{proof}
Let $\pi$ be a nonsymplectic Lefschetz fibration on $X=\product$. By Lemma \ref{nsltb},
$\pi$ is relatively
minimal, has no critical points and the fibers are tori. It is
a nontrivial bundle since otherwise it would be symplectic. If the base
space $B$ is a torus, then $X$ admits a symplectic Lefschetz fibration by Corollary
\ref{t2ovt2slf}. If $B=S^2$, then $X$ doesn't admit a symplectic structure by Lemma
\ref{t2overs2} and in particular, $M$ doesn't fiber over $S^1$ since such a fibration
would induce a symplectic form on $X$. Finally, if $genus (B) \geq 2$ and $M$ has no
fake 3--cells, then Lemma \ref{sl2dnfib} implies that $M$ doesn't fiber over $S^1$.
\end{proof}

\section{Seiberg--Witten invariants of symplectic manifolds and
     $S^1$--bundles over surfaces}

In this section we use Seiberg--Witten theory of symplectic manifolds and
$S^1$--bundles over closed, oriented surfaces to prove the following theorem which in
turn implies that the existence of a symplectic form and a Lefschetz fibration on
$\product$ is possible only if there is a symplectic Lefschetz fibration on $\product$
(Theorem \ref{cjtrlefs}). Statement (L) of Theorem \ref{mainthm} is a consequence of
this.

\begin{thm} \label{swnotlet}
Let $M$ be the total space of an oriented $S^1$--bundle over a Riemann surface $B$.
Then $X= \product$ admits a symplectic structure if and only if the bundle is trivial
or $B$ is a torus.
\end{thm}

The following theorem follows from the work of Mrowka, Ozsv\'{a}th and Yu on the SW
invariants of Seifert fibered spaces \cite{moy}. See \cite{bal} for a different (and
more elementary) approach.

\begin{thm} \label{swseif}
Let $M$ be the $S^1$--bundle over a Riemann surface $B$ of genus $g \geq 1$ with Euler
class $n \lambda $, where $\lambda $ is the (positive) generator of $H^2 (B ; \tam )$.
If $n \neq 0$, then all basic classes of $M$ are in $ \{ k \cdot \pi^* (\lambda) \ | \
0 \leq k \leq | n | -1 \} $, where $\pi$ is the bundle projection. Moreover, we have
\begin{equation} \label{sws1bun}
SW_M ( k \cdot \pi^* (\lambda) ) = \sum_{s \equiv k \ (mod \ n)} SW_{S^1 \times B} (s
\cdot pr_2^* (\lambda))\ , 
\end{equation}
where $pr_2$ is the projection $S^1 \times B \to B$.
\end{thm}

It is well-known that the Seiberg--Witten invariants of $S^1 \times B$ are given by
$$ SW_{S^1 \times B} (t) = (t -t^{-1})^{2g-2} \ ,$$
where $g$ is the genus of $B$ and the coefficient of $t^p$ on the right hand side
corresponds to the Seiberg--Witten invariant of the $Spin_c$ structure with
determinant line bundle $L$ with $c_1(L) = p \cdot pr_2^*(\lambda)$. Therefore the sum
of all Seiberg--Witten invariants of $S^1 \times B$ is $0$ if $g >1$. This sum is
preserved under twisting of the $S^1$--bundle as can be seen from (\ref{sws1bun}).

\begin{cor} \label{sumsw0} Let $M$ be as in the previous theorem and $g>1$.
Then
$$ \sum_{\alpha} SW_M(\alpha) = 0 \ ,$$
where the sum is over all $Spin_c$ structures on $M$.
\end{cor}

The following is also well-known and relates the Seiberg--Witten invariants of
$\product$ with those of $M$. For a proof see \cite{okotel}.

\begin{thm} \label{sw3-4} If $M$ is a closed, oriented 3--manifold, then
$$ SW_M (\alpha) = SW_{\product} (pr_2^* (\alpha)) $$
for any $\alpha \in H^2 (M ; \tam) $, where $pr_2$ is the projection $\product \to M$.
Moreover, if $b_+(\product) = b_1(M) > 1$, 
then all basic classes of $\product$ are pull-backs of basic classes of $M$.
\end{thm}

\proof[Proof of Theorem \ref{swnotlet}] If the bundle is trivial then $X=
T^2 \times B$ and there is a symplectic form on $X$ which is simply the sum of
symplectic forms on $T^2$ and $B$.

If $B$ is a torus, then $X$ is a torus bundle over a torus and by Theorem
\ref{t2overt2} it admits a symplectic structure.

If the bundle is nontrivial and $B$ is a sphere, then $X$ is a nontrivial
$T^2$--bundle over $S^2$ and cannot be symplectic as we proved in Lemma
\ref{t2overs2}.

From now on we will assume that the bundle is nontrivial and the genus of B is at
least 2.

By Corollary \ref{sumsw0} and Theorem \ref{sw3-4} (as $b_1(M) \geq 2 b_1(B) \geq 4$)
\begin{equation} \label{sum0}
\sum_{\alpha} SW_M (\alpha) = \sum_{\beta} SW_X(\beta) = 0\ , 
\end{equation}
where sums are over all $Spin_c$ structures on $M$ and $X$ respectively.

Assume that $X$ admits a symplectic form $\omega$. First of all, by the conditions on
equality in Theorem \ref{swsymp}, the canonical class $K= c_1(X,\omega)$ cannot be a
nonzero torsion class. On the other hand, Theorem \ref{sw3-4} and the first part of
Theorem \ref{swseif} imply that all basic classes of $X$ are torsion. Therefore the
only basic class of $X$ is $K=0$ and $SW_X(0)= \pm 1$, in particular,
$$ \sum_{\beta} SW_X(\beta) = \pm 1 \ ,$$
where the sum is over all $Spin_c$ structures on $X$. This contradicts (\ref{sum0})
hence $X$ does not admit a symplectic structure. \endproof

\begin{thm} \label{cjtrlefs} Let $M$ be a closed, oriented
3--manifold such that $\product$ admits a Lefschetz fibration and a symplectic form.
Then $\product$ admits a symplectic Lefschetz fibration or $M$ has a fake 3--cell.
\end{thm}

\begin{proof}
Let $X=\product$ admit a Lefschetz fibration and a symplectic form. Assume that there
is no symplectic Lefschetz fibration on it. Then by Lemma \ref{nsltb} it admits a
torus bundle over a Riemann surface $B$. Any such bundle should be nontrivial since
otherwise it would be symplectic. By Theorem \ref{oldmainthm}, $B$ is not a torus, and it
cannot be a sphere by Lemma \ref{t2overs2}. So $genus (B) \geq 2$. If $M$ has no fake
3--cells, then as we have seen in the proof of Lemma \ref{sl2dnfib}, $M$ is the total
space of an $S^1$--bundle over $B$ and this contradicts Theorem \ref{swnotlet}.
\end{proof}

This theorem (together with Theorem \ref{slfimpfib}) finishes the proof of statement
(L) of Theorem \ref{mainthm}.

\begin{rem}
Symplectic Lefschetz fibrations on product 4--manifolds were classified in
\cite{chemat}. As a result of our discussion, we see that nonsymplectic Lefschetz
fibrations on nonsymplectic $\product$ are nontrivial torus bundles over spherical or
hyperbolic surfaces. On the other hand, nonsymplectic Lefschetz fibrations on a
symplectic $\product$ are torus bundles over tori and by Proposition \ref{oldmainthm}
any such manifold admits a symplectic Lefschetz fibration.
\end{rem}

\section{Complex structures and Seifert fibrations on the product four--manifolds}

In this section, we use the classification of complex surfaces to prove statements (K)
and (C) of Theorem \ref{mainthm}. To prove the latter, we also use an interesting
result in Seiberg--Witten theory of complex surfaces due to Biquard. Then we consider
Seifert fibered product 4--manifolds and prove that those which admit symplectic
structures also admit either K\"{a}hler structures or torus bundles over tori. This
observation finishes the proof of Theorem \ref{mainthm}.

At this point we know exactly when the existence of a Lefschetz fibration on
$\product$ is sufficient for $M$ to fiber over the circle. Since our motivation
is to determine whether the existence of a symplectic structure on $\product$ is
sufficient for $M$ to fiber over the circle, it is quite natural to ask which
symplectic (product) 4--manifolds admit Lefschetz fibrations. This question doesn't
seem to be any easier than Conjecture T itself even though Donaldson proved that every
symplectic 4--manifold admits a Lefschetz pencil. In fact, statement (L) of Theorem
\ref{mainthm} implies that they are equivalent when $M$ has no fake 3--cells. On the
other hand, allowing multiple fibers and considering Seifert fibrations, one can still
get interesting results on Conjecture T. Seifert fibered product 4--manifolds turn out
to be closely related to complex surfaces and this is the main reason of our
discussion on complex structures on product 4--manifolds.

Now suppose that $\product$ is a closed complex surface. Since it is a spin
4--manifold its intersection form is even, so there is no exceptional sphere to
blow-down, thus it is a minimal complex surface. We are going to use the
Enriques--Kodaira classification of compact complex surfaces (see \cite{gomsti} or
\cite{barpet}) to prove the following theorem.

\begin{thm}[cf.\ Theorem 4.1 in \cite{geigo}] \label{symcom}
Let $\product$ be a closed 4--manifold.

If $\product$ admits a complex structure, then it is either an elliptic surface or of
$Class \ \rm{VII}_0$.

If $\product$ is also symplectic, then the only possibilities are the following:

{\rm(i)}\qua $\product \cong S^2 \times T^2$.

{\rm(ii)}\qua $\product$  admits a $T^2$--bundle over $T^2$.

{\rm(iii)}\qua $\product$ admits a Seifert fibration over a hyperbolic orbifold.
\end{thm}

\begin{proof}
Let $\kappa (X)$ be the Kodaira dimension of $X= \product$ as a complex surface.

{\bf{Case 1: $\kappa (X) = - \infty$}}\qua In this case $X$ is either $\kar P^2$ or
geometrically
ruled or of $Class \ \rm{VII}_0$. The complex projective plane $\kar P^2$ is
simply-connected, but $X$ is not. If
$X$ is a complex surface of $Class \ \rm{VII}_0$, then
$0=b_1(X)-1 = b_+ (X) $ hence
it cannot be symplectic. If it is geometrically ruled, then it
is the total space of a
$\kar P^1$--bundle over a Riemann surface $B$ and $ 0 = \chi (X) = \chi (\kar P^1)
\cdot \chi (B)$, hence $B$ is a torus. Moreover, $X$ is diffeomorphic to $S^2 \times
T^2$ since the total space of the nontrivial $S^2$--bundle over $T^2$ is not spin.

{\bf{Case 2: $\kappa (X) = 0$}}\qua Any minimal complex surface of Kodaira dimension 0 is
a $K3$ surface, an Enriques surface, a primary Kodaira surface, a secondary Kodaira
surface, a hyperelliptic surface or a complex torus. Since $b_1(X) \geq 1$ $X$ cannot
be a $K3$ or an Enriques surface. In three of the other four cases, $X$ is
diffeomorphic to the total space of a $T^2$--bundle over $T^2$. When $X$ is a
secondary Kodaira surface it admits an elliptic fibration over $\kar P^1$ (without
singular fibers) and $b_1(X)=1$. So in this case, $X$ cannot be symplectic because
$b_+(X)=b_1(X)-1=0$.

{\bf{Case 3: $\kappa (X) = 1$}}\qua In this case $X$ is a (properly) elliptic surface. An
elliptic fibration
on $X$ cannot have singular fibers but only
multiple fibers since the Euler characteristic of $X$ vanishes.
In particular, $X$ is a Seifert 4--manifold. While investigating
geometric structures on elliptic surfaces Wall (see \cite{wall} or \cite{wall2})
proves that the base orbifold of such a fibration must be hyperbolic if $\kappa (X) =
1$.

These are the only possibilities since every minimal surface of general type has positive Euler
characteristic, but$\chi (X) = 0$. 
\end{proof}

\begin{rem}
By a well-known result of Bogomolov \cite{tel} a complex surface of $Class$
$\rm{VII}_0$ with vanishing second Betti number is either a Hopf surface or an Inoue
surface. Since the center of the fundamental group of an Inoue surface is trivial (cf.\
Proposition 4.2 in \cite{geigo}) no Inoue surface is a product. On the other hand,
Kato's work on Hopf surfaces \cite{kato} implies that if a Hopf surface is
diffeomorphic to a product, then it must be elliptic. In particular, it is Seifert
fibered since vanishing of the Euler characteristic implies that an elliptic fibration
on a product can have no singular fibers (but only multiple ones).
\end{rem}

Recall that a closed complex surface is K\"{a}hler if and only if its first Betti
number is even. Therefore statement (K) of Theorem \ref{mainthm} is a consequence of
the following theorem.

\begin{thm} \label{kahmfib}
Let $\product$ be a closed, connected complex surface. If $b_1(M)$ is odd and $M$ has
no fake 3--cells, then $M$ is a Seifert fibered space which fibers over $S^1$.
\end{thm}

\begin{proof}
Since $b_1(X)= b_1(M) +1$ is even, $X=\product$ admits a K\"{a}hler structure. By
Theorem \ref{symcom}, $X$ is diffeomorphic to $S^2 \times T^2$ or admits a
$T^2$--bundle over $T^2$ or a properly elliptic fibration without any singular
(possibly with multiple) fibers.

If $X$ is diffeomorphic to $S^2 \times T^2$, then $M$ fibers over $S^1$ by Theorem
\ref{slfimpfib}. Moreover, the diffeomorphism between $\product$ and $S^1 \times (S^2
\times S^1)$ gives a homotopy equivalence between $M$ and $S^2 \times S^1$ and as they
both fiber over $S^1$ this homotopy equivalence must be a homeomorphism, in
particular, $M$ is a Seifert fibered space.

If $X$ admits a $T^2$--bundle over $T^2$, then by Corollary \ref{t2ovt2slf} and
Theorem \ref{slfimpfib} $M$ fibers over $S^1$ with fiber a torus and in particular it
is geometric. On the other hand, by Theorem 3 in \cite{gei} the geometric type of $M$
is $\Euc^3$, where $\Euc^n$ is $\Real^n$ with its standard metric. This implies that
$M= T^3$ (see p.446 in \cite{sco}). In particular, $M$ is Seifert fibered.

If $X$ admits a Seifert fibration over a hyperbolic orbifold $B$, then it is geometric
and the geometric type of it must be $\Euc^2 \times \Hyp^2$ by Theorem 4.5 in
\cite{wall2} as $X$ admits a K\"{a}hler structure, where $\Hyp^2$ is the hyperbolic
plane. It should be noted that there is a mistake in \cite{wall2} which was later
corrected by Kotschick in \cite{kot}; since it concerns manifolds with nonvanishing
Euler characteristic, it doesn't effect our discussion on product 4--manifolds. On the
other hand, we get the following exact sequence from the Seifert fibration
$$ 1 \To \pi_1(F) \To \pi_1(\product) \To \pi^{orb}_1 (B) \To  1 \ ,$$
where $F$ is a regular fiber and $\pi_1^{orb}(B)$ denotes the fundamental group of $B$
as an orbifold. This exact sequence leads to another one
$$ 1 \To \tam \To \pi_1(M) \To \pi^{orb}_1 (B) \To  1 \ ,$$
just as in the proof of Lemma \ref{sl2dnfib}, since $B$ is hyperbolic and its orbifold
fundamental group doesn't contain an infinite cyclic normal subgroup. So there exists
an infinite cyclic normal subgroup in $\pi_1(M)$ and $M$ is a Seifert 3--manifold by
Corollary 12.8 in \cite{hem}. (Note that as $b_1(M)$ is odd it is nonzero and $M$ is
sufficiently large.) In particular, $M$ is geometric. Since $\product$ is type $\Euc^2
\times \Hyp^2$, $M$ must be type $\Euc^1 \times \Hyp^2$, in other words the rational
Euler class of a Seifert fibration on $M$ is $0$. A generalization of Lemma
\ref{lemsl2dnfib} (e.g.\ Theorem 8.1 in \cite{neuray}) implies that $M$ fibers over
$S^1$.
\end{proof}

In order to prove statement (C) of Theorem \ref{mainthm} we use the following result
of Biquard (cf.\ Th\'{e}or\`{e}me 8.2 in \cite{biq}):

\begin{thm} \label{biqthm}
A properly elliptic non--K\"ahler surface   admits no symplectic structure.
\end{thm}

\proof[Proof of Statement (C) in Theorem \ref{mainthm}] We have seen in
Theorem \ref{symcom} that if $X=\product$ admits a complex and a symplectic structure,
then there are three possibilities. The product $S^2 \times T^2$ admits a K\"{a}hler structure
hence if $X=S^2 \times T^2$, then $M$ fibers over $S^1$ by Theorem \ref{kahmfib}. If
$X$ admits a $T^2$--bundle over $T^2$, then $M$ fibers over $S^1$ by Corollary
\ref{t2ovt2slf} and Theorem \ref{slfimpfib}.  If $X$ is a properly elliptic surface,
then it has to be K\"{a}hler by Theorem \ref{biqthm} hence $M$ fibers over $S^1$ by
Theorem \ref{kahmfib}. \endproof

The following is a well-known theorem. For a nice proof see \cite{woo}.
\begin{thm} \label{seicom}
If $M$ is a closed, oriented Seifert fibered space, then $\product$ admits a complex
structure.
\end{thm}

\begin{prop}\label{sympsei}
Let $M$ be a closed, oriented 3--manifold with no fake 3--cells. Suppose $\product$
admits a symplectic structure and a Seifert fibration. Then $\product$ admits a
K\"{a}hler structure or a $T^2$--bundle over $T^2$.
\end{prop}

\begin{proof}
We have the following short exact sequence coming from the Seifert fibration
$$ 1 \To \pi_1(F) \To \pi_1(\product) \xrightarrow{\pi_{\#}} \pi^{orb}_1 (B)
\To  1 \ , $$ 
where $F$ is a generic fiber, $ \pi_1^{orb} (B)$ denotes the fundamental
group of $B$ as an orbifold and $\pi$ is the projection map of the fibration. Let $u$
be a generator of $\pi_1(S^1 \times \{ pt \} ) $ in $\pi_1 (\product)$ as in the proof
of Lemma \ref{sl2dnfib}.

First assume that $\pi_{\#} (u)$ is nontrivial in $\pi_1^{orb}(B)$. Then it generates
an infinite, cyclic, normal subgroup (cf.\ proof of Lemma \ref{sl2dnfib}). Existence of
such a subgroup in $\pi_1^{orb}(B)$ is possible only if $B$ is a nonsingular orbifold
of genus 1, i.e. a torus. So the Seifert fibration we have is in fact a $T^2$--bundle
over $T^2$.

Now assume $u \in ker (\pi_{\#} )$. Then as in the proof of Theorem \ref{kahmfib} we
have
$$ 1 \To \tam \To \pi_1(M) \To \pi^{orb}_1 (B) \To  1 \ .$$
In particular, there is an infinite cyclic normal subgroup of $\G$. Since $X$ admits a
symplectic structure $b_+ (X) \geq 1$ and so is $b_1(M)$. This implies that $M$ is
sufficiently large. Therefore we can use Corollary 12.8 in \cite{hem} to conclude that
$M$ is a Seifert fibered space. So $\product$ admits a complex structure by Theorem
\ref{seicom}, hence it admits a K\"{a}hler structure or a $T^2$--bundle over $T^2$ as
in the proof of statement (C).
\end{proof}

This proposition (together with Theorem \ref{kahmfib} and Corollary \ref{t2ovt2slf})
finishes the proof of Theorem \ref{mainthm}.

\section{Geometry of $M$ and structures on $\product$}

During the course of our proof of Theorem \ref{mainthm} we made observations on the
interaction between various structures and fibrations on $M$ and $\product$. In this
section, we recall some of those observations and use them to prove a couple of
theorems on the relation between the geometry of $M$ and $\product$.

Throughout this section we will assume that $M$ is a closed, connected and oriented
3--manifold with no fake 3--cells.

In the proof of Proposition \ref{sympsei} we used the existence of a symplectic
structure on $\product$ to conclude that $b_+(\product)=b_1(M) > 0$. Note that $b_1(M)
> 0$ implies that $M$ is sufficiently large.

\begin{thm}
If $\product$ is Seifert fibered and $M$ is sufficiently large, then $M$ admits a
nonhyperbolic geometric structure.
\end{thm}

\begin{proof}
As in the proof of Proposition \ref{sympsei} we look at the homotopy sequence of the
Seifert fibration. There are two different cases depending on the image of a generator
$u$ of $\pi_1(S^1 \times \{ pt \} ) \subset \pi_1(\product )$:

If $u$ is in the kernel, then we have an infinite cyclic normal subgroup in
$\pi_1(M)$. Since $M$ is sufficiently large, Corollary 12.8 in \cite{hem} implies that
$M$ is a Seifert fibered space.

If $u$ is not in the kernel, then $\product$ admits a $T^2$--bundle over $T^2$, in
particular it is symplectic. Hence (e.g.\ by (L) of Theorem \ref{mainthm}) $M$ fibers
over the circle with fiber a torus. By Theorem 5.5 in \cite{sco} $M$ is geometric of
type $\Euc^3$, $Nil^3$ or $Sol^3$.

It is now clear that in any case $M$ is geometric but not hyperbolic.
\end{proof}

As we mentioned before if $M$ is Seifert fibered, then $\product$ admits a complex
structure. If $M$ is geometric of type $Sol^3$, then $\product$ is obviously geometric
of type $\Euc^1 \times Sol^3$ and as a consequence $\product$ doesn't admit any
complex structure \cite{wall}.

On the other hand, Theorem \ref{symcom} says that if $\product$ admits a complex
structure, then it is either of $Class \ \rm{VII}_0$ or an elliptic surface and in any
case, by the remark following Theorem \ref{symcom} $\product$ is Seifert fibered.

This discussion leads us to the following conclusion which is a partial converse of
the well-known Theorem \ref{seicom}.

\begin{thm}
If $\product$ admits a complex structure and $M$ is sufficiently large, then $M$ is a
Seifert fibered space.
\end{thm}

\Addressesr

\end{document}

%% file: agtout.tex
\def\ifplaintex{\expandafter\ifx\csname documentclass\endcsname\relax}

\def\gtp{{\mathsurround=0pt\it $\cal G\mskip-2mu$eometry \&\ 
$\cal T\!\!$opology $\cal P\!$ublications}}  

\def\Addressesr{\bigskip
{\small \parskip 0pt \leftskip 0pt \rightskip 0pt plus 1fil \def\\{\par}
\sl\theaddress\par
\medskip
\rm Email:\stdspace\tt\theemail\hfill\rm Received:\qua\receiveddate \par}}

\def\recd{{\small Received:\qua\receiveddate\ifx\reviseddate\relax
\else\qquad Revised:\qua\reviseddate\fi\par}} 

\def\lognumber#1{\def\thelognumber{#1}}
\def\volumenumber#1{\def\thevolumenumber{#1}}
\def\volumeyear#1{\def\thevolumeyear{#1}}
\def\papernumber#1{\def\thepapernumber{#1}}
\def\pagenumbers#1#2{\def\startpage{#1}\def\finishpage{#2}}
\def\published#1{\def\publishdate{#1}}

\def\received#1{\def\receiveddate{#1}}

\def\accepted#1{\def\accepteddate{#1}}
\def\asciititle#1{\def\theasciititle{#1}}
\def\covertitle#1{\def\thecovertitle{#1}}
\def\asciiauthors#1{\def\theasciiauthors{#1}}

\def\coverauthors#1{\def\thecoverauthors{#1}}
\long\def\asciiabstract#1{\long\def\theasciiabstract{#1}}
\def\asciikeywords#1{\def\theasciikeywords{#1}}

\let\\\par\let\thelognumber\relax\let\thevolumenumber\relax
\let\thepapernumber\relax\let\thevolumeyear\relax\let\startpage\relax
\let\finishpage\relax\let\publishdate\relax\let\receiveddate\relax
\let\reviseddate\relax\let\accepteddate\relax\let\theasciititle\relax
\let\thecovertitle\relax\let\theasciiauthors\relax
\let\theasciiabstract\relax\let\theasciikeywords\relax

\let\thecoverauthors\relax\let\theasciiemail\relax

\ifplaintex
\font\logobig=cmssbx10 scaled 3836
\font\logomed=cmssbx10 scaled 2557
\else
\font\logobig=cmssbx10 scaled 4200
\font\logomed=cmssbx10 scaled 2800
\fi

\long\def\makeagttitle{   
\count0=\startpage
\agt\hfill      
\hbox to 45truept{\vbox to 0pt{\vglue -13truept{\logomed A\kern -.37em{\logobig 
T}\kern -.38em G}\vss}\hss}
\break
{\small Volume \thevolumenumber\ (\thevolumeyear)
\startpage--\finishpage\nl
Published: \publishdate}

\vglue .25truein

{\parskip=0pt\leftskip 0pt plus
1fil\def\\{\par\smallskip}{\Large\bf\thetitle}\par\medskip} \vglue
0.05truein

{\parskip=0pt\leftskip 0pt plus 1fil\def\\{\par}{\sc\theauthors}
\par\medskip}%
 
\vglue 0.03truein 

{\small\leftskip 25truept\rightskip 25truept{\bf Abstract}\stdspace\theabstract

{\bf AMS Classification}\stdspace\theprimaryclass
\ifx\thesecondaryclass\relax\else; \thesecondaryclass\fi\par
{\bf Keywords}\stdspace \thekeywords\par}\vglue 7truept

}   

\ifplaintex
\font\phead=cmsl9 scaled 950
\font\pnum=cmbx10 scaled 913
\font\pfoot=cmsl9 scaled 950
\headline{\vbox to 0pt{\vskip -4.5mm\line{\small\phead\ifnum
\count0=\startpage ISSN 1472-2739 (on-line) 1472-2747 (printed)
\hfill {\pnum\folio}\else\ifodd\count0\def\\{ }%
\ifx\theshorttitle\relax\thetitle\else\theshorttitle\fi\hfill{\pnum\folio}
\else\def\\{ and }{\pnum\folio}\hfill\ifx\theshortauthors\relax\theauthors
\else\theshortauthors\fi\fi\fi}\vss}}
\footline{\vbox to 0pt{\vglue 0mm\line{\small\pfoot\ifnum\count0=\startpage
\copyright\ \gtp\hfill\else
\agt, Volume \thevolumenumber\ (\thevolumeyear)\hfill\fi}\vss}}
\else
\font=cmsl9 scaled 1050
\font\lnum=cmbx10 
\font=cmsl9 scaled 1050
\makeatletter
\def\@oddhead{{\small\lhead\ifnum\count0=\startpage ISSN 1472-2739 
(on-line) 1472-2747 (printed)\hfill {\lnum\number\count0}\else\ifodd\count0
\def\\{ }\ifx\theshorttitle\relax \thetitle \else\theshorttitle\fi\hfill
{\lnum\number\count0}\else\def\\{ and }{\lnum\number\count0}
\hfill\ifx\theshortauthors\relax 
\theauthors\else\theshortauthors\fi\fi\fi}}\def\@evenhead{\@oddhead}
\def\@oddfoot{\small\lfoot\ifnum\count0=\startpage\copyright\ \gtp\hfill\else
\agt, Volume \thevolumenumber\ (\thevolumeyear)\hfill\fi}
\def\@evenfoot{\@oddfoot}
\makeatother
\fi
\let\maketitlepage\makeagttitle
\let\makeshorttitle\maketitlepage
\let\maketitle\maketitlepage

\newwrite\gtoutfile
\long\gdef\makeheadfile{  
{\def\\{, }\def\s{ }
\immediate\openout\gtoutfile head.xxx
\immediate\write\gtoutfile{To: math@arxiv.org}
\immediate\write\gtoutfile{Subject: put OR rep NNNNN:ppppp}
\immediate\write\gtoutfile{--text follows this line--}
\immediate\write\gtoutfile{Proxy-for: \ifx\theasciiauthors\relax
\theauthors\else\theasciiauthors\fi\s<\ifx\theasciiemail\relax\theemail\else\theasciiemail\fi>}
\immediate\write\gtoutfile{\noexpand\\}
\immediate\write\gtoutfile{Authors: \ifx\theasciiauthors\relax
\theauthors\else\theasciiauthors\fi}
{\def\\{ }\immediate\write\gtoutfile{Title: \ifx\theasciititle\relax
\thetitle\else\theasciititle\fi}}
\immediate\write\gtoutfile{Subj-class: GT or SG, GR etc}
\immediate\write\gtoutfile{MSC-class: \theprimaryclass\ifx\thesecondaryclass\relax\else, \thesecondaryclass\fi}
\immediate\write\gtoutfile{Journal-ref: Algebr. Geom. Topol. \thevolumenumber\s
(\thevolumeyear) \startpage-\finishpage}
\immediate\write\gtoutfile{Comments: Published by Algebraic and
Geometric Topology at}
\immediate\write\gtoutfile{\s\s\s  http://www.maths.warwick.ac.uk/agt/AGTVol\thevolumenumber/agt-\thevolumenumber-\thepapernumber.abs.html}
\immediate\write\gtoutfile{\noexpand\\}
\immediate\write\gtoutfile{}
\ifx\theasciiabstract\relax
\immediate\write\gtoutfile{\theabstract}\else
\immediate\write\gtoutfile{\theasciiabstract}\fi
\immediate\write\gtoutfile{}
\immediate\write\gtoutfile{\noexpand\\}
\immediate\write\gtoutfile{}
\immediate\closeout\gtoutfile}}  

\def\maketitlepage{\makeagttitle\makeheadfile}
\let\makeshorttitle\maketitlepage
\let\maketitle\maketitlepage

\def\ifplaintex{\expandafter\ifx\csname documentclass\endcsname\relax}

\def\gtp{{\mathsurround=0pt\it $\cal G\mskip-2mu$eometry \&\ 
$\cal T\!\!$opology $\cal P\!$ublications}}  

\def\Addressesr{\bigskip
{\small \parskip 0pt \leftskip 0pt \rightskip 0pt plus 1fil \def\\{\par}
\sl\theaddress\par
\medskip
\rm Email:\stdspace\tt\theemail\hfill\rm Received:\qua\receiveddate \par}}

\def\recd{{\small Received:\qua\receiveddate\ifx\reviseddate\relax
\else\qquad Revised:\qua\reviseddate\fi\par}} 

\def\lognumber#1{\def\thelognumber{#1}}
\def\volumenumber#1{\def\thevolumenumber{#1}}
\def\volumeyear#1{\def\thevolumeyear{#1}}
\def\papernumber#1{\def\thepapernumber{#1}}
\def\pagenumbers#1#2{\def\startpage{#1}\def\finishpage{#2}}
\def\published#1{\def\publishdate{#1}}

\def\received#1{\def\receiveddate{#1}}

\def\accepted#1{\def\accepteddate{#1}}
\def\asciititle#1{\def\theasciititle{#1}}
\def\covertitle#1{\def\thecovertitle{#1}}
\def\asciiauthors#1{\def\theasciiauthors{#1}}

\def\coverauthors#1{\def\thecoverauthors{#1}}
\long\def\asciiabstract#1{\long\def\theasciiabstract{#1}}
\def\asciikeywords#1{\def\theasciikeywords{#1}}

\let\\\par\let\thelognumber\relax\let\thevolumenumber\relax
\let\thepapernumber\relax\let\thevolumeyear\relax\let\startpage\relax
\let\finishpage\relax\let\publishdate\relax\let\receiveddate\relax
\let\reviseddate\relax\let\accepteddate\relax\let\theasciititle\relax
\let\thecovertitle\relax\let\theasciiauthors\relax
\let\theasciiabstract\relax\let\theasciikeywords\relax

\let\thecoverauthors\relax\let\theasciiemail\relax

\ifplaintex
\font\logobig=cmssbx10 scaled 3836
\font\logomed=cmssbx10 scaled 2557
\else
\font\logobig=cmssbx10 scaled 4200
\font\logomed=cmssbx10 scaled 2800
\fi

\long\def\makeagttitle{   
\count0=\startpage
\agt\hfill      
\hbox to 45truept{\vbox to 0pt{\vglue -13truept{\logomed A\kern -.37em{\logobig 
T}\kern -.38em G}\vss}\hss}
\break
{\small Volume \thevolumenumber\ (\thevolumeyear)
\startpage--\finishpage\nl
Published: \publishdate}

\vglue .25truein

{\parskip=0pt\leftskip 0pt plus
1fil\def\\{\par\smallskip}{\Large\bf\thetitle}\par\medskip} \vglue
0.05truein

{\parskip=0pt\leftskip 0pt plus 1fil\def\\{\par}{\sc\theauthors}
\par\medskip}%
 
\vglue 0.03truein 

{\small\leftskip 25truept\rightskip 25truept{\bf Abstract}\stdspace\theabstract

{\bf AMS Classification}\stdspace\theprimaryclass
\ifx\thesecondaryclass\relax\else; \thesecondaryclass\fi\par
{\bf Keywords}\stdspace \thekeywords\par}\vglue 7truept

}   

\ifplaintex
\font\phead=cmsl9 scaled 950
\font\pnum=cmbx10 scaled 913
\font\pfoot=cmsl9 scaled 950
\headline{\vbox to 0pt{\vskip -4.5mm\line{\small\phead\ifnum
\count0=\startpage ISSN 1472-2739 (on-line) 1472-2747 (printed)
\hfill {\pnum\folio}\else\ifodd\count0\def\\{ }%
\ifx\theshorttitle\relax\thetitle\else\theshorttitle\fi\hfill{\pnum\folio}
\else\def\\{ and }{\pnum\folio}\hfill\ifx\theshortauthors\relax\theauthors
\else\theshortauthors\fi\fi\fi}\vss}}
\footline{\vbox to 0pt{\vglue 0mm\line{\small\pfoot\ifnum\count0=\startpage
\copyright\ \gtp\hfill\else
\agt, Volume \thevolumenumber\ (\thevolumeyear)\hfill\fi}\vss}}
\else
\font=cmsl9 scaled 1050
\font\lnum=cmbx10 
\font=cmsl9 scaled 1050
\makeatletter
\def\@oddhead{{\small\lhead\ifnum\count0=\startpage ISSN 1472-2739 
(on-line) 1472-2747 (printed)\hfill {\lnum\number\count0}\else\ifodd\count0
\def\\{ }\ifx\theshorttitle\relax \thetitle \else\theshorttitle\fi\hfill
{\lnum\number\count0}\else\def\\{ and }{\lnum\number\count0}
\hfill\ifx\theshortauthors\relax 
\theauthors\else\theshortauthors\fi\fi\fi}}\def\@evenhead{\@oddhead}
\def\@oddfoot{\small\lfoot\ifnum\count0=\startpage\copyright\ \gtp\hfill\else
\agt, Volume \thevolumenumber\ (\thevolumeyear)\hfill\fi}
\def\@evenfoot{\@oddfoot}
\makeatother
\fi
\let\maketitlepage\makeagttitle
\let\makeshorttitle\maketitlepage
\let\maketitle\maketitlepage

\newwrite\gtoutfile
\long\gdef\makeheadfile{  
{\def\\{, }\def\s{ }
\immediate\openout\gtoutfile head.xxx
\immediate\write\gtoutfile{To: math@arxiv.org}
\immediate\write\gtoutfile{Subject: put OR rep NNNNN:ppppp}
\immediate\write\gtoutfile{--text follows this line--}
\immediate\write\gtoutfile{Proxy-for: \ifx\theasciiauthors\relax
\theauthors\else\theasciiauthors\fi\s<\ifx\theasciiemail\relax\theemail\else\theasciiemail\fi>}
\immediate\write\gtoutfile{\noexpand\\}
\immediate\write\gtoutfile{Authors: \ifx\theasciiauthors\relax
\theauthors\else\theasciiauthors\fi}
{\def\\{ }\immediate\write\gtoutfile{Title: \ifx\theasciititle\relax
\thetitle\else\theasciititle\fi}}
\immediate\write\gtoutfile{Subj-class: GT or SG, GR etc}
\immediate\write\gtoutfile{MSC-class: \theprimaryclass\ifx\thesecondaryclass\relax\else, \thesecondaryclass\fi}
\immediate\write\gtoutfile{Journal-ref: Algebr. Geom. Topol. \thevolumenumber\s
(\thevolumeyear) \startpage-\finishpage}
\immediate\write\gtoutfile{Comments: Published by Algebraic and
Geometric Topology at}
\immediate\write\gtoutfile{\s\s\s  http://www.maths.warwick.ac.uk/agt/AGTVol\thevolumenumber/agt-\thevolumenumber-\thepapernumber.abs.html}
\immediate\write\gtoutfile{\noexpand\\}
\immediate\write\gtoutfile{}
\ifx\theasciiabstract\relax
\immediate\write\gtoutfile{\theabstract}\else
\immediate\write\gtoutfile{\theasciiabstract}\fi
\immediate\write\gtoutfile{}
\immediate\write\gtoutfile{\noexpand\\}
\immediate\write\gtoutfile{}
\immediate\closeout\gtoutfile}}  

\def\maketitlepage{\makeagttitle\makeheadfile}
\let\makeshorttitle\maketitlepage
\let\maketitle\maketitlepage